\documentclass[12pt,a4paper]{article} 
\usepackage{amsmath,pstricks,latexsym,theorem,epsfig} 
\newtheorem{Theorem}{Theorem}[section] 
\newtheorem{Definition}{Definition}[section] 
\newtheorem{Proposition}{Proposition}[section] 
\newtheorem{Lemma}{Lemma}[section] 
\newtheorem{Corollary}{Corollary}[section] 
 
\theorembodyfont{\rmfamily}

\def\OO{\mathcal{O}}
\def \R{I\!\!R}

\def\11{1\!\!1}

\newcommand{\ba}{\begin{array}} 
\newcommand{\ea}{\end{array}}

\begin{document} 
\title{\bf Partial Regularity for Stationary Solutions to    Liouville-Type Equation in dimension 3} 
\author{{Francesca Da Lio\footnote{ Dipartimento di Matematica Pura ed Applicata, Universit\`a degli Studi di Padova. Via Trieste 63, 35121,Padova, Italy, e-mail: dalio@mah.unipd.it. }} \footnote
{Department of Mathematics, ETH Z\"urich, R\"amistrasse 101, 8092 Z\"urich, Switzerland.}}
%\date{} 

\maketitle 
 \begin{abstract}
In dimension $n=3$, we prove that the singular set of any stationary solution to the Liouville equation $-\Delta u=e^u$,  which belongs to $W^{1,2}$, has Hausdorff dimension at most $1$.
 \end{abstract}
\section{Introduction} 

The regularity theory for nonlinear elliptic equations has a long history. It is beyond the scope of the present work to describe even part of it and we better refer the reader chapter 14 in \cite{Tay}  for a presentation of this theory and further references.\par  Typical examples of nonlinear elliptic problems under study are the  semilinear elliptic problems of the form
  \begin{equation}\label{semil}
 Lu=f(x,u,\ldots,\nabla^{m-1}u)\,,
 \end{equation}
where the function $u$ is defined  on some open subset of $\R^n$,  $L$ is a linear elliptic operator of order $m$ and where the nonlinear operator $f $ involves derivatives of $u$ up to  order $m-1\,.$\par
Once we fixe the dimension $n$ of the underlying space and the function space ${\cal{V}}$, to which the solution $u$ is assumed to belong, equations like (\ref{semil}) can be classified in three categories :  the sub-critical, critical and super-critical equations. \par

These three categories of equations (which depend on the choice of $n$ and $\cal V$) are characterized as follows.  Starting from the fact that  $u$ belongs to ${\cal{V}}$, one can estimate the nonlinear part $f(x, u,\ldots,\nabla^{m-1}u)$, if in addition $u$ is a solution to (\ref{semil}), this implies that $Lu$ belongs to some function space ${\cal{W}}$ (which is usually larger than the space ${\cal{V}}$ itself)\,. Sub-critical (respectively super-critical) equations,  are the one for which the information $Lu \in {\cal{W}}$  implies, through elliptic regularity theory, that $u$ belongs to a function space  which is strictly smaller (respectively strictly larger) and which has different homogeneities than the original ${\cal{V}}$\,. In turn, critical equations are the one which are neither sub-critical nor supercritical in the above sense \,.  \par

It is well known that solutions to subcritical equations of the form (\ref{semil}) for a smooth $f$ and smooth $L$  are in fact smooth. This is a consequence of the standard {\it bootstrap argument}\,. In contrast with the subcritical situation, solutions to a given critical equation either can all be proven to be smooth or can have non trivial singular sets (that is, non removable singularities). These results then depend on the nature of the nonlinearity $f$. 

For example,  in dimension $n=2$, when ${\cal V}=W^{1,2}(B^2,{\R})$, the equation
  \begin{equation}
  \label{a1}
  -\Delta u=|\nabla u|^2
  \end{equation}
is critical. Indeed, plugging the information $u\in W^{1,2}(B^2,{\R})$ into $f(\nabla u)=|\nabla u|^2$, one obtains that $\Delta u\in L^1$ which itself implies that $\nabla u$ is in$L^{2,\infty}$,  the weak$-L^2$ space, which has the same homogeneity as $L^2$. Thus, in a some  sense, we are back to  the initial situation and this shows that the equation is critical. Observe that this critical equation, when $n=2$ and ${\cal V}=W^{1,2}(B^2,{\R})$, admits singular solutions such as $\log\log\frac{1}{r}$.
  
In contrast to the above situation, one can consider the equation
  \begin{equation}
  \label{a2}
  -\Delta u=u_x\wedge u_y
  \end{equation}
which is again critical in dimension $n=2$ when ${\cal V}=W^{1,2}(B^2,{\R}^3)$, but this time any solution can be shown to be smooth (see for instance \cite{Hel}).\par
  
Finally, in dimension $n = 3$ and when ${\cal V}=W^{1,2}(B^3,{\R}^3)$, this equation is  super-critical and the existence result of T. Rivi\`ere \cite{Riv} of everywhere discontinuous harmonic maps in $W^{1,2}(B^3,S^2)$ has annihilated  all hope of having a partial regularity result for solution to this super-critical semilinear equation.  \par
     
When the equation has a variational structure, namely when the equation is the Euler-Lagrange equation of a functional, it makes sense to restrict our attention to the subspace of solutions which are {\it stationary}. That is,  one considers the critical points to the functional which are also critical with respect to perturbations of the domain (see Definition \ref{statdef} below and see also \cite{Hel}).  A consequence of this {\it stationarity assumption} is that the solution satisfies an identity (which is in fact a conservation law) which, in the most studied cases, can be converted into a {\em monotonicity formula}. In most of the cases which have been studied so far, this monotonicity formula implies that the solution $u$ belongs to some Morrey type space ${\cal M}$, which is much smaller than the original space ${\cal{V}}$. In the good cases, replacing  ${\cal{V}}$ by ${\cal M}$ makes the problem critical and this allows one to obtain a partial regularity result for the stationary solutions (see for instance \cite{Hel,E} for harmonic maps, and \cite{P} when the nonlinearity is $u^\alpha$ with $\alpha$ greater than the critical exponent )\,.\par
  
The aim of this paper is to present an alternative approach to the partial regularity theory when the stationary assumption cannot be converted in a {\it monotonicity formula}. We illustrate this method by applying it to the famous Liouville equation in dimension $n=3$ 
\begin{equation}\label{emden}
-\Delta u=e^{u}~~\mbox{in $\Omega$\,.}
\end{equation} \par

Throughout the paper, $\Omega\subseteq\R^3\,$ denotes an open set, $u$ is a scalar function and $\nabla u$ and $\nabla^2 u$ denote respectively the gradient and Hessian matrix of $u$.  In dimension $n=2$, the geometric meaning of equation (\ref{emden})  is well known and it corresponds to the problem of finding metrics $g$, which are conformally equivalent to the flat metric, and which have constant Gauss curvature.
In dimension  $n\ge 3$, equation (\ref{emden}) arises in the modeling of several physical phenomena such as the theory of isothermal gas sphere and   gas combustion \,. 
\par
A function $u$ is said to be a weak solution of (\ref{emden}) in $\Omega$ if for all $\varphi\in C_c^\infty(\Omega)$ it satisfies
\begin{equation}\label{weak}
-\int_{\Omega} u\Delta \varphi dx=\int_{\Omega} e^u \varphi dx\,.
\end{equation}
We now recall the definition of stationary solution.\par
\begin{Definition}\label{statdef}
A weak solution of (\ref{emden}) is said to be stationary if  it satisfies
\begin{equation}\label{deren} \frac{d}{dt} E(u(x+t X))_{|_{t=0}}=0\,,\end{equation}
for all smooth vector fields $X$ with compact support in $\Omega$, where
$$E(u)=\frac{1}{2}\int_{\Omega}|\nabla u|^2dx -\int_{\Omega} e^u dx\,.$$
\end{Definition}
Computing (\ref{deren}) for weak solutions in $W^{1,2}(\Omega)$ we find that for any smooth vector field $X$ the following identity holds
\begin{equation}\label{stationary}
\int_{\Omega} \left[ \frac{\partial u}{\partial x_i}\frac{\partial u}{\partial x_k}\frac{\partial X^k}{\partial x_i}-\frac{1}{2}|\nabla u|^2\frac{\partial
X^i}{\partial x_i}+e^u \frac{\partial X^i}{\partial x_i} \right] dx=0\,,
\end{equation}
This identity can be also understood as a conservation law (see again \cite{Hel}).\par

Arguing  as in  \cite{HL}, we insert in  (\ref{stationary}) the vector field
   $X^\delta=x\varphi^\delta(|x|)$ where 
   $$
   \left\{\begin{array}{ll}
   \varphi^\delta(s)=1 &\mbox{
    if $ s<r$}\\
     \varphi^\delta(s)=1+\frac{r-s}{r\delta} & \mbox{if $r\le s\le r+r\delta$}\\
     \varphi^\delta(s)=0 &\mbox{
    if $ s> r$}\,.
     \end{array}
    \right.$$
 \par  
  After some calculations,  we let $\delta\to 0$ and deduce that, for almost every $r>0$, the following  formula holds
  \begin{equation}\label{monotonicity}
   \frac{1}{r}\int_{ B_{r}}(\frac{1}{2}|\nabla u|^2-3e^u )dx= \frac{1}{2}\int_{\partial B_{r}} |\nabla_T u|^2 dx-\frac{1}{2}\int_{\partial B_{r}}
  |\frac{\partial u}{\partial r}|^2 dx-
  \int_{\partial B_{r}} e^u dx\,.
   \end{equation}
This can also be written as follows
   \begin{equation}
   \frac{d}{dr}\left[\frac{1}{r}\int_{B_r}(|\nabla u|^2-6 e^u)dx \right]=\frac{2}{r}\int_{\partial B_{r}}
  |\frac{\partial u}{\partial r}|^2 dx-\frac{4}{r}\int_{\partial B_{r}}  e^udx\,.
  \end{equation}
 Unlike the cases of stationary solutions to super-critical semilinear equations which have mainly be considered so far, the formula (\ref{monotonicity}) does not seem to provide any monotonicity information,  any uniform bound neither for the term 
 $\frac{1}{r}\int_{ B_{r}(x_0)}|\nabla u|^2 dx$ nor for  $\frac{1}{r}\int_{ B_{r}(x_0)}e^{u} dx\,.$ As already mentioned, the main contribution of the present work is to present an alternative approach to the partial regularity theory in abscence of Monotonicity and Morrey type estimates. Our approach is inspired by the technique introduced by Fang-Hua Lin and Tristan Rivi\`ere in \cite{LR1} in the context of Ginzburg-Landau equations. This technique based on some kind of dimension reduction argument. More precisely,  applying Fubini's Theorem one first extracts  ``good" 2 dimensional slices
to get  estimates of the some suitable quantities, then one restricts these quantities to these slices (whose dimension is such that the non-linearity $e^u$ becomes critical for $W^{1,2}$) and obtain some estimates in interpolation spaces : the Lorentz spaces $L^{2,\infty}-L^{2,1}$.  Finally, the stationarity condition (\ref{monotonicity}) can be used to "propagate" these estimates from the slices (basically the boundary of balls) into the domain bounded by the slices (the balls themseves). \par
 
Now we state our main result.
 
  \begin{Theorem}\label{regularity}
  Assume that $u\in W^{1,2}(\Omega)$ is a stationary solution of (\ref{emden}), such that $e^u\in L^1(\Omega)$. Then there exists an open set $\OO\subset\Omega$ such that
  $$u\in C^{\infty}(\OO)~~ \mbox{and} ~~{\cal{H}}_{dim}(\Omega\setminus\OO)\le 1\,$$
  where ${\cal{H}}_{dim}$ denotes the dimensional Hausdorff measure\,.
  \end{Theorem} 
  
  It is an open question whether such a partial regularity result is optimal or not (the same question holds for instance also for stationary harmonic maps). What is known is that stationary solutions to (\ref{emden}) can have singularities. Indeed the function
  $u(x)=\log(\frac{2}{|x|^2})$ satisfies $-\Delta u = e^u$ but is not bounded\,.\par
  
Our approach and the above result should also hold for the more general class of equations of the form $-\Delta u = V(x) e^u$ where $V(x)$ is some
smooth given potential. For the sake of simplicity, we have chosen to focus our attention on the case where $V\equiv 1$ in order to the keep the technicalities as low as possible and make the paper more "readable".
We recall  that, in dimension $n=2$, the regularity of weak solutions  to the equation  (\ref{emden}), starting from the hypothesis that $u$ is in $W^{1,2}$ , is a straighforward consequence of the Moser-Trudinger
inequality (see \cite{GT} ). Still in dimension $n=2$, a $L^{\infty}$ estimate for solutions in  $L^{1}(\Omega)$  to  the equation  (\ref{emden}), starting from the hypothesis that  $e^u\in L^1(\Omega)\,,$ has been obtained  by Brezis \& Merle \cite{BM}. Finally,  in \cite{BLO} the authors prove some {\em a priori} estimates for solutions of (\ref{emden}) in any dimension but under the stronger assumption
  $e^u$ is in some ad-hoc Morrey Space which makes the problem critical.\par
    
\section{Preliminary Estimates of the Energy }\label{prel}
In this Section we are going to prove some preliminary estimates.\par
We first  introduce some notations and recall the definition of Hausforff measure.\par
 For $x_0\in\Omega$, $r>0$ we will denote by $B_{r,x_0}$ or simply by  $B_r$ the ball $B(x_0,r)$ centered at $x_0$ and with radius $r$. 
Given $A\subset\R^3$ 
 we denote by $|A|$ its Lebegue measure and by ${\cal{H}}^s(A)$ its $s$   dimensional  Hausdorff measure\,.
   \par
 We   recall (see e.g \cite{EG}) the definition of the $s$- dimensional Hausforff measure ${\cal{H}}^s$ in $\R^n$, with $0\le s\le n\,.$
    For any $\delta>0$ and for any $A\subseteq\R^n$ we set
    $$
    {\cal{H}}^s_\delta(A)=\inf\{\sum_{i}w_s r^s_i~:A\subseteq\cup_{i}B_{r_i}.~r_i<\delta,~\forall i\}$$
    where $w_s=\frac{\pi^{s/2}}{\Gamma(1+s/2)}$ and the infimum is taken over all contable collections of ball $\{B_{r_i}\}$ covering the set $A$ and having radii
    $r<\delta\,.$ The $s$-dimensional  Hausdorff measure is then defined as
    $${\cal{H}}^s(A)=\lim_{\delta\to 0}{\cal{H}}^s_\delta(A)\,.$$
  Given $x_0\in\Omega$ and $0<r< d(x_0,\partial\Omega)$ 
   we introduce the following energy 
  \begin{equation}\label{energy}
  {\cal{E}}_{r,x_0}(u)=\frac{1}{r}\int_{B_r}|\nabla u|^2 dx+\frac{1}{r}\int_{B_r} e^u dx\,.
  \end{equation}
  and   set
  $$
  (u)_{x_0,r}:=\frac{1}{|\partial B_{r,x_0}|}\int_{\partial B_{r,x_0}} u(y) dy\,,
  $$
  The key result to prove Theorem \ref{regularity} is the following assertion about the energy (\ref{energy}).

  \begin{Theorem}\label{decreasing}
   There exist constants $\eta,\beta\in (0,1)$  such that for every $x_0\in\Omega$ and $0<r<d(x_0,\partial\Omega),$ 
\begin{equation}\label{small}
   {\cal{E}}_{r,x_0}(u)\le \eta ~~\mbox{and}~~ (u)_{r,x_0} {\cal{E}}_{r,x_0}(u)\le \eta
   \end{equation}
   imply
   $$
    {\cal{E}}_{s,y}(u)\le C s^{\beta}
    $$
    for all  $y $ in a neighborhood of $x_0$ for all  $s\le r$, and  for some  $C$ depending on $r\,$ and independent on $y,s\,.$
    \end{Theorem}

  In order to prove Theorem \ref{decreasing} we   need to give some definitions and to show a series of preliminary results\,.  
  \par
    
  We start with recalling the definition of the weak $L^2$ space (or Marcinkievicz space $L^{2,\infty})\,,$ (see \cite{Ste})\,.\par
  The space $L^{2,\infty}(\Omega)$ is   defined as the space of functions $f\colon \Omega\to\R$ such that
   $$\sup_{\lambda\in\R}\lambda|\{x:|f|(x)\ge\lambda\}|^{1/2}<+\infty$$
    The dual space of $L^{2,\infty}(\Omega)$ is the Lorentz space $L^{2,1}(\Omega)$ whose norm is equivalent
   $$
   ||f||_{2,1}\simeq\int_0^\infty 2|\{x: |f|(x)\ge s\}|ds\,.$$
   In dimension $2$ we have the following property\,:
   $W^{1,1}(\Omega)$   continuously embedds  in $L^{2,1}(\Omega)\,,$ (see e.g \cite{SW,Hunt,T1} )\,.\par 
  We next  recall  a  result proved by   Lin \& Riviere in \cite{LR1} in the framework of Ginzurg-Landau functionals, which will play a crucial role
  in getting  estimates of the energy (\ref{energy})\,.
  \begin{Lemma}\label{LinRiv}{(Choice of a ``good" slice)}[Lemma A.2, \cite{LR1}]
  For any $g\in L^1(B_1)$,  if we denote 
  $$f(x)=\int_{B_1}\frac{g(y)}{|y-x|^2}dy\,,$$
 then  for every $\delta>0$ there exists a subset $E_\delta\subseteq (0,1)$ and $|E_\delta|>1-\delta$ such that 
  for all $\rho\in E_\delta$ we have
 $$ ||f||_{L^{2,\infty}(\partial B_\rho)}\le C_\delta ||g||_{L^1(B_1)} \,,$$
 where $C_\delta$ depends only on $\delta$\,.
 \end{Lemma}
  
   Since $e^u,|\nabla u|^2$ are in $L^{1}(\Omega),$ we will suppose in the sequel without restriction that 
    $$\int_{\partial B_r} (e^u+|\nabla u|^2 )dx<+\infty\,.$$
   We   decompose $u-(u)_{x_0,r}$ as the sum of two functions solving to different Dirichlet Problems. More presisely we
    write $u-(u)_{x_0,r}=v+w$, with $v$ and $w$ satisfying respectively
    \begin{equation}\label{DPv}
  \left\{\begin{array}{ll}
  -\Delta v=e^u~~&\mbox{in $B_r$}\\
  v=0~~&\mbox{on $\partial B_r\,;$}
  \end{array}\right.
  \end{equation}
   \begin{equation}\label{DPw}
  \left\{\begin{array}{ll}
  -\Delta w=0~~&\mbox{in $B_r$}\\
 w=u- (u)_{x_0,r}~~&\mbox{on $\partial B_r$}
  \end{array}\right.
  \end{equation} 
  In the next two subsections we estimate  $\frac{1}{\rho r}\int_{B_{\rho r}}e^u dx$ and
 $\frac{1}{\rho r}\int_{B_{\rho r}}|\nabla u|^2 dx$ in function of the energy ${\cal{E}}_{r,x_0}(u)\,,$ by using specific properties satisfied by $v$ and $w\,.$
 
  %Estimates of $e^u$
   \subsection{Estimates of $e^u$}
   In this subsection we are going to estimate $\frac{1}{\rho r}\int_{B_{\rho r}}e^u dx$ in function of the energy ${\cal{E}}_{r,x_0}(u)\,.$ 
   More precisely we prove the following Theorem.
   \begin{Theorem}\label{estexpu}
 For all $ \alpha\in (0,1)$,  there exist constants $\eta\in (0,1)$,
  $0<\rho_1<\rho_2<1$  such that
   \begin{equation}\label{smallbis} {\cal{E}}_{r,x_0}(u)\le \eta\,,\end{equation}
   implies  
   \begin{equation}
   \frac{1}{\rho r}\int_{B_{\rho r}}e^u dx\le \alpha\,{\cal{E}}_{r,x_0}(u)\,,
   \end{equation}
   for   every $\rho\in [\rho_1,\rho_2]$ \,.
   \end{Theorem}
   {\bf Proof.} We split the proof in several steps.\par
   {\bf Step 1.} We start by estimating  $e^w\,.$
    
  \begin{Lemma}\label{expw} 
  The function $\omega=e^w$ satisfies  
  \begin{equation}\label{decrmean}
  \frac{1}{\rho r}\int_{ B_{\rho r}}e^w dx\le \rho^2\frac{1}{r}\int_{ B_r}e^w dx\,.
  \end{equation}
  \end{Lemma}
  {\bf Proof.} We observe that $e^w$ is subharmonic\par
  $$-\Delta e^w=-|\nabla w|^2 e^w \le 0\,.$$
  A well known fact of sub-harmonic functions is that their mean value on a ball is 
       a nonincreasing function with respect to radius of the ball, namely the following holds for every $\rho\in (0,1)\,$
  $$
  \frac{1}{|B_{\rho r}|}\int_{B_{\rho r}} e^w dx \le  \frac{1}{|B_{r}|}\int_{B_{r}} e^w dx\,.$$
 This clearly  implies
  $$
  \frac{1}{\rho}\int_{B_{\rho r}} e^w dx \le  \rho^2\frac{1}{r}\int_{B_{r}} e^w dx \,,
  $$
and we conclude\,. ~~\hfill$\Box$\par
 {\bf Step 2.}
 \begin{Proposition}
For all $\rho\in (0,1)$ and $x,y\in B_{\rho r}$ we have
$$|w(x)-w(y)|\le C\eta^{1/2} (2\rho r)\,,
$$
for some $C$ depending only on the dimension of the space\,.
\end{Proposition}
{\bf Proof.}
We set $\bar w_r=w(rx+x_0)-(u)_{r,x_0},$
$\bar w_r$ satisfies 
 \begin{equation} 
  \left\{\begin{array}{ll}
  -\Delta \bar w_r=0~~&\mbox{in $B_{1}$}\\
 \bar w_r=u_r-(u)_{r,x_0}~~&\mbox{on $\partial B_{1}\,.$}
  \end{array}\right.
  \end{equation} 
  Standard elliptic estimates and  Poinca\'re-Wirtinger Inequality imply that  for all  $\rho\in(0,1)$ and for foer $C>0$
 \begin{eqnarray*}
 && ||\bar w_r||_{C^1(B_{\rho})}\le C(\int_{\partial B_1}(u_r-(u)_{r,x_0})^2 dx)^{1/2}\le C||\nabla u_r||_{L^2(\partial B_{1})}\,.
\end{eqnarray*}
    Thus for every $x,y\in B_{\rho r}$ we have
  \begin{equation}
  |w(x)-w(y)|^2\le C|x-y|^2\int_{\partial B_{r}}|\nabla u(z)|^2dz\le C|x-y|^2\int_{B_{r}}|\nabla u(z)|^2dz
  \end{equation}
  From  the   assumption  (\ref{smallbis}) on the energy  it follows that 
\begin{equation}
  |w(x)-w(y)| \le C(2\rho r)\eta^{1/2},
  \end{equation}
  and we can conclude.~\hfill$\Box$\par
   
   \begin{Proposition}\label{estvs}
 The function $v$ satisfies
 \begin{equation} \label{estv}
 \left(\frac{1}{r^3}\int_{B_r} v^2(x) dx\right)^{1/2}\le C\eta \,.
 \end{equation}
 \end{Proposition}
 {\bf Proof of Proposition \ref{estvs}.} We recall that  
 $$v(x)=\int_{B_r}e^{u(y)} G(x,y) dy,$$
 where $G(x,y)$ is the Green function on the ball which satisfies 
 $$|G(x,y)|\le C\frac{1}{|x-y|} ~~\mbox{and}~~ |\nabla_x G(x,y)|\le C\frac{1}{|x-y|^2}\,,$$(see e.g. \cite{GT}).\par
 Thus
 \begin{equation} 
 \left( \int_{B_1} v^2(rx) dx\right)^{1/2}\le Cr^2\int_{ B_1}e^{u(rx)} dx\,.
 \end{equation}
 By a change of variable we get 
 $$\left(\frac{1}{r^3}\int_{B_r} v^2(x) dx\right)^{1/2}\le C\frac{1}{r}\int_{ B_r}e^u dx\,.$$
 By applying assumption  (\ref{smallbis}) we get 
 (\ref{estv}) and we conclude\,.~~\hfill $\Box$
 
 ~~\hfill$\Box$
 
  {\bf Step 3.} 
  From  Proposition \ref{estvs} it follows 
  \begin{equation}\label{esteta}
  (|\{x\in B_{\rho r}: v\geq \eta^{1/2}\}|\eta)^{1/2}\le C\eta r^{3/2} 
 \end{equation}
  and 
   \begin{equation}\label{esteta2} |\{x\in B_{\rho r}: v\geq \eta^{1/2}\}|\le C\eta^{1/2} r^3= C\eta^{1/2}\frac{1}{\rho^3}(\rho r)^3\,.\end{equation}
  Now take $\lambda>0$ (that we will   determine  later)  and we set
  \begin{eqnarray*}&& I^1_{\eta}:=\{x:~~v(x)\le \eta\}\,,\\&&I^1_{\lambda}:=\{x:~~v(x)\ge \lambda\}\,,\\&&I^2_{\lambda,\eta}=\{x:~~\eta^{1/2}\le v(x)\le \lambda\}\,.\\
  \end{eqnarray*}
  The following estimates holds.
  \begin{eqnarray*}
  \frac{1}{r}\int_{ B_r}ve^u dx&=&\frac{1}{r}\int_{ B_r}v(-\Delta v) dx\\&=&\frac{1}{r}\int_{ B_r}|\nabla v|^2\le \frac{1}{r}\int_{ B_r}|\nabla u|^2\,.
  \end{eqnarray*}
  Thus 
 \begin{eqnarray}\label{i1}
  \frac{1}{\rho r}\int_{ B_{\rho r}\cap I^1_\lambda}e^u dx&\le& \frac{1}{\rho r\lambda}\int_{ B_{\rho r}\cap I^1_\lambda}|\nabla u|^2dx\\
  &\le &\frac{1}{\rho r\lambda}\int_{ B_r} |\nabla u|^2dx\,.\nonumber
   \end{eqnarray}
We also have  
\begin{eqnarray}\label{i2}
  \frac{1}{\rho r}\int_{ B_{\rho r}\cap I^1_\eta}e^u dx&\le& e^{\eta^{1/2}}\frac{1}{\rho r }\int_{ B_{\rho r}}e^{w}dx\,.
\end{eqnarray}
By Proposition \ref{expw}, estimate (\ref{esteta2}) and the fact that $e^{w(0)}\le \frac{1}{|B_{\rho r}|}\int_{  B_{\rho r}} e^{w(x)}dx\,,$ we get
\begin{eqnarray}\label{i3}
  \frac{1}{\rho r}\int_{ B_\rho\cap I^2_{\lambda,\eta}}e^ve^w dx&\le & 
e^\lambda \frac{1}{\rho r}\int_{  B_{\rho r}\cap I^2_{\lambda,\eta}}e^w dx\le
e^\lambda e^{w(0)}\frac{1}{\rho r}\int_{ B_{\rho r}\cap I^2_{\lambda,\eta}}e^{w(y)-w(0)} dx\nonumber\\
&\le& e^\lambda  e^{C\eta^{1/2}}|\{x~ v\ge \eta^{1/2}\}|\frac{1}{(\rho r)^3}\int_{  B_{\rho r}} e^wdx \nonumber\\
&\le& Ce^\lambda   e^{C\eta^{1/2}}\eta^{1/2}\frac{r}{\rho^2}\frac{1}{\rho r}\int_{  B_{\rho r}} e^wdx\,.
\end{eqnarray}
 By combining the above estimates (\ref{i1}),  (\ref{i2}) and  (\ref{i3})  we finally  get
\begin{eqnarray*}
\frac{1}{\rho r}\int_{ B_{\rho r}} e^u dx&=&
\frac{1}{\rho r}\int_{B_{\rho r}\cap I^1_\eta} e^u dx +\frac{1}{\rho r}\int_{B_{\rho r}\cap I^1_\lambda} e^u dx\\
&&~~+\frac{1}{\rho r}\int_{ B_{\rho r}\cap I^2_{\lambda,\eta}} e^u dx\\
&\le&e^{\eta^{1/2}}\frac{1}{\rho r}\int_{  B_{\rho r}} e^w dx+\frac{1}{\rho r\lambda}\int_{ B_{\rho r}}|\nabla u|^2dx\\
&&~~+Ce^{\lambda} e^{C\eta^{1/2}}\eta^{1/2}\frac{r}{\rho^2}\frac{1}{\rho r}\int_{ B_{\rho r}} e^wdx\\
&\le& e^{\eta^{1/2}}\rho^2\frac{1}{r}\int_{  B_{ r}} e^w dx +
\frac{C}{\rho \lambda}\frac{1}{r}\int_{ B_{  r}}|\nabla u|^2dx\\
&&~~+Ce^{\lambda} e^{C\eta^{1/2}}\eta^{1/2}r\frac{1}{r}\int_{ B_r} e^wdx\\
&=&[C e^\lambda e^{C\eta^{1/2}}\eta^{1/2}r+\rho^2 e^{\eta^{1/2}}]\frac{1}{r}\int_{ B_r} e^wdx\\&&~~+ 
 \frac{1}{\rho \lambda}\frac{1}{r}\int_{ B_{r}}|\nabla u|^2dx
\end{eqnarray*}
 Now we   fix the interval $[\rho_1,\rho_2]$ where we make $\rho$ vary, and  constants  $ \lambda,\eta\,.$
 \par
 We  consider any  $0<\alpha<<1$.  We first choose $\rho$ such that
 $$ e\rho ^2<\frac{\alpha}{3}\,.$$
 Thus we    take $ \rho_1,\rho_2 $ satsfying $0<\rho_1<\rho_2<\frac{\sqrt{\alpha}}{\sqrt{3 e}}\,.$
  Then  we choose 
  $\lambda$ large enough so that
$$
\frac{1}{\rho_1 \lambda}<\frac{\alpha}{3}
$$
and finally  we choose $\eta$ small enough so that
$$  C e^{C\eta^{1/2}}e^\lambda\eta^{1/2}<\frac{\alpha}{3}$$
We observe that $\int_{ B_r} e^wdx\le \int_{ B_r} e^udx$ , being $v$ nonnegative by the Maximum Principle. Thus with these choices of the constants $\rho_1,\rho_2$,
 $\eta$ and $\lambda$ 
we obtain
$$
\frac{1}{\rho r}\int_{ B_{\rho r}} e^u dx\le
\alpha[\frac{1}{r}\int_{ B_r} e^u dx+\frac{1}{r}\int_{ B_r} |\nabla u|^2 dx]\,,
$$
 for all $\rho\in[\rho_1,\rho_2]\,.$ Thus we can    conclude.~\hfill$\Box$
 
  \subsection{Estimate of $\nabla u$}
  
  %Estimate gradient u
  In   this subsection we are going to estimate $\frac{1}{r\rho}\int_{B_{\rho r}}|\nabla u|^2 dx$ in function of
  ${\cal{E}}_{r,x_0}(u)$, 
  for $\rho\in[\rho_1,\rho_2]\subseteq(0,1)$   $\rho_1,\rho_2$ being the constants determined in Theorem \ref{estexpu}\,. 
   \begin{Theorem}[Estimate gradient of $u$]\label{Thgradu}
      For almost every   $\rho\in[\rho_1,\rho_2]$  we have
    \begin{eqnarray}\label{decrgradu}
  \frac{1}{\rho r}\int_{B_{\rho r}}|\nabla u|^2 dx&\le&
    C({\cal{E}}_{r,x_0}(u)) \left({\cal{E}}_{r,x_0}(u)  + (u)^{+}_{r,x_0}{\cal{E}}_{r,x_0}(u)\right)+\gamma ({\cal{E}}_{r,x_0}(u))\,.\nonumber\\
    \end{eqnarray}
    for some $C>0$ and  $0<\gamma<1$ independent on $r\,.$
  \end{Theorem}
  %%%%%%%%%%%%%%%%%%%%%%%%%%%%%%%%%%%%%%%%%%%%%%%%%%%%%%
  {\bf   Proof.}
  We split the proof in several steps. 
  
  %%%%%%%%%%%%%%%%%%%%%%%%%%%%%%%%%%%%%%%%%%%%%%%%%
 
  {\bf Step 1.}
       We start  by showing that $\nabla v$ and $\nabla w$ are orthogonal in $B_r$\,.
    \begin{Lemma}\label{orth}
  The following estimate holds
  $$\ \int_{B_{r}} |\nabla u|^2 dx= \int_{B_{r}} |\nabla v|^2 dx+ \int_{B_{r}} |\nabla w|^2 dx\,.
  $$
  \end{Lemma}
  {\bf Proof.} 
  Let $\nu(x)$ denote the exterior normal versor to $\partial B_r$ at the point $x\in\partial B_r$.
  We have
 $$\int_{B_{r}} \nabla u\cdot\nabla v dx=  \int_{\partial B_{r}} v\nabla w\cdot\nu- \int_{B_{r}} v\Delta w=0\,.
 $$
~~\hfill$\Box$\par
  {\bf Step 2.}    Estimate of   $\nabla w .$
   \begin{Proposition}\label{gradw}
For every $\rho\in (0,1)$ we have    
   \begin{equation}\label{estgradw}
   \frac{1}{r\rho}\int_{B_{\rho r}} |\nabla w|^2 dx\le \rho^2\frac{1}{r}\int _{B_r} |\nabla w|^2 dx\,.
   \end{equation}
   \end{Proposition}
   {\bf Proof.} We observe that $w\in C^\infty(B_r)$ and $\omega=|\nabla w|^2$ satisfies
   $$
   -\Delta\omega(x)\le 0,~~\mbox{in $B_r$\,.}$$
   The conclusion follows as Lemma \ref{expw} \,.
    ~\hfill$\Box$

    {\bf Step 3.}    Estimate of   $ \nabla v \,.$
    We start by showing   some intermediate estimates.\par

   \begin{Proposition}\label{lapl}
   For some $C>0$ (independent on $r$) we have
   \begin{eqnarray}\label{logcond}
  \int_{B_r}|\Delta v|\log(2+|\Delta v|) dx&\le & C  (\int_{B_r} (e^{u(x)} +|\nabla u(x)|^2dx  + (u)^{+}_{r,x_0} \int_{B_r} e^{u(x)}dx)\,.
  \end{eqnarray}
  \end{Proposition}
  {\bf Proof of Proposition \ref{lapl}.} 
  We have \par
  \begin{eqnarray}\label{estlogcond}
   \int_{B_r}|\Delta v|\log(2+|\Delta v|) dx& = &  \int_{B_r}e^u\log(2+e^u) dx\nonumber\\&\le & 
    C \int_{B_r}e^u(1+u^+) dx\\ &\le &   C [\int_{B_r}e^u dx+(u)^{+}_{r,x_0} \int_{B_r}e^u dx\nonumber \\
    &+& \int_{B_r}e^u v^+ dx+\int_{B_r}e^u w^+ dx]\nonumber\end{eqnarray}
   We estimate the two last terms of (\ref{estlogcond}).\par
   \begin{eqnarray}\label{vterm}
   \int_{B_r}e^u v^+ dx&=&\int_{B_r}-\Delta v v^+ dx=\int_{B_r}|\nabla v|^2 dx\nonumber\\
   &\le & \int_{B_r}|\nabla u|^2 dx\,.
   \end{eqnarray}
   \begin{eqnarray}\label{wterm}
   \int_{B_r}e^u w^+ dx&=&  \int_{B_r}-\Delta v w^+ dx \nonumber \\
   &=&-\int_{B_r}div(\nabla v w^+) dx+\int_{B_r}\nabla v\cdot\nabla w^+\nonumber\\
   &=&-\int_{\partial B_r}\frac{\partial v}{\partial \nu} w^+ dx\nonumber \\&\le& 2 \int_{\partial B_r}|\nabla v|^2 dx 
   + 2\int_{\partial B_r}[(u-(u)_{r,x_0})^+]^2 dx\nonumber \\
   &\le& 2 \int_{\partial B_r} |\nabla v|^2 dx+2 \int_{\partial B_r} |\nabla u|^2 dx\,.\nonumber\\
   &\le & 4 \int_{  B_r} |\nabla u|^2 dx
   \end{eqnarray}
   In the estimate (\ref{wterm}) we use the fact that $v\in H^2(\partial B_r)$ , being $u,w\in W^{1,2}(\partial B_r)\,$ and thus
   $$
   \int_{\partial B_r} |\nabla v|^2 dx\le \int_{ B_r} |\nabla v|^2 dx\,.$$
   
    By combining (\ref{estlogcond}), (\ref{vterm}), (\ref{wterm})  we get
     \begin{eqnarray*} 
  \int_{B_r}|\Delta v|\log(2+|\Delta v|) dx&\le &C [\int_{B_r}e^u dx+(u)^{+}_{r,x_0} \int_{B_r}e^{u(x)} dx\\
  &+& 5\int_{\partial B_r} |\nabla u|^2 dx]\\ &\le &
  C  (\int_{B_r} \left(e^{u(x)} +|\nabla u(x)|^2dx  + (u)^{+}_{r,x_0} \int_{B_r} e^{u(x)}dx\right)\,.
  \end{eqnarray*}
     Thus we can conclude.~\hfill$\Box$
     \begin{Corollary}\label{hessv} 
     We have $\nabla^2 v\in L^1(B_r)$ and
     $$
     ||\nabla^2 v||_{ L^1(B_r)}\le C  \left(\int_{B_r} (e^{u(x)} +|\nabla u(x)|^2dx  + (u)^{+}_{r,x_0} \int_{B_r} e^{u(x)}dx\right)\,.$$
     \end{Corollary}
     {\bf Proof.}   Calderon-Zygmund theory (see e.g. \cite{SW}) yields   that if
     $\int_{B_r}|\Delta v|\log(2+|\Delta v|) dx<+\infty$ then $\nabla^2 v\in L^1(B_r)$ and
     $$
     ||\nabla^2 v||_{ L^1(B_r)}\le C\int_{B_r}|\Delta v|\log(2+|\Delta v|) dx\,.$$
     Thus the result  follows directly from Proposition \ref{lapl} and we conclude\,.~\hfill $\Box$
     
We can now use  Lemma \ref{LinRiv} to prove  the following result.   
      \begin{Proposition}\label{corLinRiv}
      For   every $\delta>0$ small enough, there
    exists a subset $E_\delta\subseteq [\rho_1,\rho_2]$ and $|E_\delta|>\rho_2-\rho_1-\delta$ such that 
  for all $\rho\in E_\delta$ we have
 $$ ||\nabla v||_{L^{2,\infty}(\partial B_{\rho r})}\le C_\delta  \int_{B_r} e^{u(x)}  dx\,,$$
 where $C_\delta$ depends only on $\delta$\,.
  \end{Proposition}
  {\bf Proof of Proposition \ref{corLinRiv}.} As we observe in Proposition \ref{estvs},   we can write 
  $$v(x)=  \int_{B_r} e^{u(y)}G(x,y)dy\, ,$$
   \par
   where $G(x,y)$ is the Green function on $B_r\,.$ Since $|\nabla_x G(x,y)|\le C\frac{1}{|x-y|^2} $ we have
    $$|\nabla v(x)|\le  C\int_{B_r} \frac{e^{u(y)}}{|x-y|^2} dy\, .$$
   Lemma \ref{LinRiv} yields   that  for  every $\delta>0$ there
    exists a subset $E_\delta\subseteq [\rho_1,\rho_2]$ and $|E_\delta|>\rho_2-\rho_1-\delta$ such that 
  for all $\rho\in E_\delta$ we have
  $$||\nabla v ||_{L^{2,\infty}(\partial B_{\rho r})}\le C_\delta \int_{B_r}e^{u(x)} dx\,,$$
  and we conclude\,.~~\hfill\mbox{$\Box$}

    \begin{Proposition}\label{estgradv}
    For every $\delta>0$ there
    exists a subset $E_\delta\subseteq [\rho_1,\rho_2]$ and $|E_\delta|>\rho_2-\rho_1-\delta$ such that 
  for almost every  $\rho\in E_\delta$ we have
    $$
    \int_{\partial B_{\rho r}} |\nabla v|^2 dx\le
     C\left( {\cal{E}}_{r,x_0}(u)\right)\left( {\cal{E}}_{r,x_0}(u)  + (u)^{+}_{r,x_0}{\cal{E}}_{r,x_0}(u)\right)\,,
    $$
    with $C$ depending on $\delta$ and the dimension\,.
    \end{Proposition}
    {\bf Proof.}
    Since $\nabla^2 v\in L^1(B_r)$ by Fubini Theorem for almost every  $\rho\in [0,1]$ we have
     $\nabla^2 v\in L^1(\partial B_{\rho r})$. By the embedding  of  the space 
    $W^{1,1}(  \partial B_{\rho r})$ into
      $L^{2,1}( \partial B_{\rho r})$ we have $\nabla v\in L^{2,1}(\partial B_{\rho r})$ as well  and the following estimate holds\par
    \begin{eqnarray*}
    ||\nabla v||_{L^{2,1}(\partial B_{\rho r})}&\le& C ||\nabla^2 v||_{L^1(\partial B_{\rho r})} \\
    &\le& \frac{C}{r}||\nabla^2 v||_{L^1(  B_r)}\le \frac{C}{r}(\int_{B_r}|\nabla^2 v|dx  )\\
    &\le&  C ({\cal{E}}_{r,x_0}(u)  + (u)^{+}_{r,x_0} {\cal{E}}_{r,x_0}(u) )\,.
    \end{eqnarray*}
    Now by using the duality  between  $L^{2,\infty}$ and $L^{2,1} $ and Proposition  \ref{corLinRiv} , we get
    \begin{eqnarray*}
    \int_{\partial B_{\rho r}} |\nabla v|^2 dx&\le& ||\nabla v||_{L^{2,\infty}(\partial B_{\rho r})}||\nabla v||_{L^{2,1}(\partial B_{\rho r})}
    \\
    &\le& C 
    (\int_{B_r}(e^u +|\nabla u|^2 dx))( {\cal{E}}_{r,x_0}(u)  + (u)^{+}_{r,x_0}\frac{1}{r}\int_{B_r} e^{u(x)} dx)\\
    &\le& C {\cal{E}}_{r,x_0}(u)( {\cal{E}}_{r,x_0}(u)  + (u)^{+}_{r,x_0}{\cal{E}}_{r,x_0}(u)) \,.
    ~~~\mbox{$\Box$}
    \end{eqnarray*}
     
  {\bf Step 4.}
    From Proposition \ref{gradw} and Fubini Theorem it follows that for almost every $\rho\in(0,1)$ 
  $$
  \int_{\partial B_{\rho r}}|\nabla w|^2 dx\le \frac{\rho^2}{r} \int_{  B_{ r}}|\nabla w|^2 dx\,.$$
  Thus  for almost every $\rho\in E_\delta$ ($E_\delta$ is as in Proposition \ref{estgradv}) we have
  \begin{eqnarray*}
  \int_{\partial B_{\rho r}}|\nabla u|^2 dx&\le& 2\int_{\partial B_{\rho r}}|\nabla w|^2 dx+2\int_{\partial B_{\rho r}}|\nabla v|^2 dx\\
  &\le&2\frac{\rho^2}{r} \int_{  B_{ r}}|\nabla w|^2 dx\\
  &&~~+C  ( {\cal{E}}_{r,x_0}(u))( {\cal{E}}_{r,x_0}(u)  + (u)^{+}_{r,x_0}{\cal{E}}_{r,x_0}(u))\,.
  \end{eqnarray*}
   By applying  formula (\ref{monotonicity}) to $u$ in the ball $B_{\rho r}$ and  Theorem \ref{estexpu}
  we obtain the following estimate for almost every $\rho\in E_\delta\subseteq[\rho_1,\rho_2]$
  \begin{eqnarray*} 
 \frac{1}{\rho r}\int_{B_{\rho r}} |\nabla u|^2 dx&\le &  \int_{\partial B_{\rho r}} |\nabla u|^2 dx+\frac{6}{\rho r}\int_{B_{\rho r}}e^u dx \\
 &\le& 2\rho^2 {\cal{E}}_{r,x_0}(u)+ 2C  {\cal{E}}_{r,x_0}(u)({\cal{E}}_{r,x_0}(u)+(u)^{+}_{r,x_0}{\cal{E}}_{r,x_0}(u))\\
  &&~~+\frac{6}{\rho r}\int_{B_{\rho r}}e^u dx\\
  &\le& 2\rho^2 {\cal{E}}_{r,x_0}(u)+ 2C  {\cal{E}}_{r,x_0}(u)({\cal{E}}_{r,x_0}(u)+(u)^{+}_{r,x_0}{\cal{E}}_{r,x_0}(u))\\
  &&~~+6\alpha {\cal{E}}_{r,x_0}(u)\\
  &\le& (2\rho^2 +6\alpha){\cal{E}}_{r,x_0}(u)+2C  {\cal{E}}_{r,x_0}(u)({\cal{E}}_{r,x_0}(u)+(u)^{+}_{r,x_0}{\cal{E}}_{r,x_0}(u))
   \end{eqnarray*}
   where $\alpha$ is the constant appearing in Theorem \ref{estexpu}\,.
   We remark that we can always choose $\rho_2$ and $\alpha$ in Theorem \ref{estexpu} in  such a way that $2\rho^2 +6\alpha<\gamma<1\,.$
  Thus we can conclude the proof of Theorem \ref{Thgradu}\,.\hfill$\Box$
  
    \section{Proofs of Theorem \ref{decreasing} and Theorem \ref{regularity}}
     
 In this Section we give the proof of Theorem \ref{decreasing} and  Theorem \ref{regularity}\,.
 We start by giving an estimate of the mean value $(u)_{r,x_0}$.
 \begin{Lemma}
 For all $0<r<s\le 1$ the following estimate holds 
 \begin{eqnarray}\label{estmean}
 (u)_{r,x_0}&\le & (u)_{s,x_0}+\frac{1}{r}\int_{B_r} e^{u(x)}dx\nonumber\\
 &&~~- \int_{B_s\setminus B_r} \frac{e^{u(x)}}{|x-x_0|}dx\,.
 \end{eqnarray}
 \end{Lemma}
 {\bf Proof.} 
 One can   check that in the sense of distribution the following estimate holds.
 \begin{equation}\label{dermean}
 \frac{d}{dr}  (u)_{r,x_0}=-\frac{1}{r^2}\int_{B_r} e^{u(x)}dx<0\,
 \end{equation}
Integrating (\ref{dermean}) between $r$ and $s$ we get 
 \begin{eqnarray}\label{estmeaneq}
 (u)_{r,x_0}&= & (u)_{s,x_0}+\frac{1}{r}\int_{B_r} e^{u(x)}dx-\frac{1}{s}\int_{B_s} e^{u(x)}dx\nonumber\\
 &&~~+\int_{B_s\setminus B_r} \frac{e^{u(x)}}{|x-x_0|}dx\nonumber\,,\\
 \end{eqnarray}
 and we conclude\,.~~\hfill$\Box$par
 \medskip
 {\bf Proof of Theorem \ref{decreasing}.}
 We split the proof in several steps.
 \par
 {\bf Step 1.}  
 By combining     Theorem \ref{estexpu} and Theorem \ref{Thgradu}   we can find $\rho\in[\rho_1,\rho_2]$ (independent on $r$) such that
 \begin{eqnarray}\label{energyest1}
 {\cal{E}}_{\rho r,x_0}(u)&\le& \gamma {\cal{E}}_{r,x_0}(u)\\
 &&~~+C {\cal{E}}_{r,x_0}(u) \left({\cal{E}}_{r,x_0}(u)  + (u)^{+}_{r,x_0}{\cal{E}}_{r,x_0}(u)\right)\,.
  \end{eqnarray}
 Indeed  we observe that up to choosing $\eta,\gamma$ and $\rho_1,\rho_2$   smaller, the constant  $\rho_1$ always satisfies (\ref{energyest1})\,. \par
 We set $\tau_j=\rho^j r$, $a_j={\cal{E}}_{\tau_j,x_0}(u)$ and $u_j=(u)_{\tau_j,x_0}\,.$
 First of all we have
 \begin{eqnarray}\label{estmean2}
 u_j&\le& u_0+\frac{1}{\tau_0}\int_{B_{\tau_0}} e^u dx\\
 &&~~+\frac{1}{\tau}\sum_{k=0}^{j}\frac{1}{\tau^k}\int_{B_{\tau_j}} e^u dx\nonumber\\
 &\le& u_0+\frac{2}{\tau}\sum_{k=0}^{j}{\cal{E}}_{\tau_k,x_0}(u)\nonumber
 \end{eqnarray}
 By plugging (\ref{estmean2}) in (\ref{energyest1}) we get
 \begin{eqnarray}\label{energyest2}
 a_{j+1}&\le& \gamma a_j +C a_j\left(a_j+a_ju_0+\frac{2}{\rho_1 r} a_j(\sum_0^j a_k)\right) 
 \end{eqnarray}
 The recursive formula (\ref{energyest2}) implies that if $\eta$ is small enough then  
 $$ a_j\le a_0{\bar \beta}^j\,,$$ 
 for some
 $0<\bar \beta<1$
 We deduce that  for all $0\le s\le r$ we have
 $${\cal{E}}_{s,x_0}(u)\le C s^{ \beta}\,,$$
with $0< \beta<1$ and $C$ is a positive constant that may depend on $r\,.$
 \par
 
 {\bf Step 2.} {\bf Claim:} The maps
 $$ \Omega\times(0,1]\to\R~~(x_0,r)\mapsto {\cal{E}}_{r,x_0}(u) $$
 and
 $$\Omega\times(0,1]\to\R,~(x_0,r)\mapsto (u)_{r,x_0}
 $$
 are continuous.
 \par
 {\bf Proof of the claim.}
 The continuity of $(x_0,r)\mapsto {\cal{E}}_{r,x_0}(u)$ follows from the fact that 
   $e^u$ and $|\nabla u|^2$ are in $L^1(\Omega)$\,.  \par
    
   The continuity of the map $x_0\mapsto (u)_{r,x_0}$ follows from  the fact it can be represented as the composition 
  of  the following three continuous maps.\par
  The first map is ~: 
        $\Omega\to W^{1,2}(B_{r,0})$, $x_0\mapsto u(x-x_0)$; \par
      The second map is the trace operator: $$W^{1,2}(B_{r,0})\to W^{1/2,2}(\partial B_{r,0}),~~
  u\mapsto u|_{\partial B_{r,0}}\,.$$\par \par 
  The third one is the  bounded linear  operator $$ H^{1/2}(\partial B_{r,0})\to \R,~~
  w\mapsto\frac{1}{|\partial B_{r,0}|}\int_{\partial B_{r,0}} w dx\,.$$ 
  \par
  Finally  we use the fact that for $0<r<1$ we have
  \begin{eqnarray}\label{estmeanbis}
 (u)_{r,x_0}&= & (u)_{1,x_0}+\frac{1}{r}\int_{B_r} e^{u(x)}dx- \int_{B_1} e^{u(x)}dx\nonumber\\
 &&~~+\int_{B_1\setminus B_r} \frac{e^{u(x)}}{|x-x_0|}dx \,.
 \end{eqnarray}
  Since the right hand side of (\ref{estmeanbis}) is continuous with respect to $(r,x_0)\in (0,1]\times\Omega$, we conclude that
  $(r,x_0)\mapsto (u)_{r,x_0}$ is continuous as well.\par
  
  {\bf Step 4.}
  By the continuity of $(x,r)\mapsto {\cal{E}}_{r,x}(u)$ and $(x,r)\mapsto (u)_{r,x}$  we can conclude that
  up to the choice of a smaller $\eta,$ we can find $\varepsilon>0$   such that for all
  $y\in B_{\varepsilon,x_0}$ and for  $s\in(r-\varepsilon,r+\varepsilon)$ we have
  $${\cal{E}}_{s,y}(u)<\eta ,~~\mbox{and} ~~(u)_{s,y}{\cal{E}}_{s,y}(u)< \eta\,.$$
  Finally  by Theorem \ref{estexpu}, Theorem \ref{Thgradu}  and Step 1 we get
  $${\cal{E}}_{s,y}(u)\le C s^{ \beta}\,,$$
  for all
  $y\in B_{\varepsilon,x_0}$ and for $s\le s_0$, $s_0$ being a constant independent on $y$ (actually by changing $C$ we could  choose $s_0=r$). Thus we can conclude the proof of Theorem \ref{decreasing}\,. ~~\hfill $\Box$
  \par\medskip
  {\bf Proof of Theorem \ref{regularity}.}
  Set
  \begin{eqnarray*}
  \OO&=&\{x\in\Omega~: {\cal{E}}_{r,x}(u)<\eta ~\mbox{and}~ (u)_{r,x}\,{\cal{E}}_{r,x}(u)<\eta\\
  &&~~~~~\mbox{for some $0<r<dist(x,\Omega)$}\}\,.
  \end{eqnarray*}
  From Theorem \ref{decreasing} it follows that $\OO$ is open. Moreover $u\in C^{\beta/2}(\OO)$, (see e.g. Giaquinta \cite{Gia}),
   and routine elliptic regularity theory then
  proves that $u\in C^{\infty}(\OO)\,.$  \par
  We set
  $$A_1=\{x\in\Omega~: {\cal{E}}_{r,x}(u)\ge \eta ~~\mbox{for all $0<r<dist(x,\Omega)$}\}$$
  and
  $$A_2=\{x\in\Omega~: (u)_{r,x}{\cal{E}}_{r,x}(u)\ge\eta ~~\mbox{for all $0<r<dist(x,\Omega)$}\}\,.
  $$
  We have 
   $$V=\OO^c=A_1\cup A_2\,.$$
   Next we show that ${\cal{H}}^1(A_1)=0$ and ${\cal{H}}^{1+\alpha}(A_2)=0\,$ for any $\alpha>0\,.$\par
   {\bf1.} ${\cal{H}}^1(A_1)=0$:~~let $x\in A_1$. By definition we have
   \begin{equation}\label{cond}
    {\cal{E}}_{r,x}(u)\ge\eta 
   \end{equation}
    for all $0<r<dist(x,\Omega)$\,. Now fix $\delta>0$ and set 
    \begin{eqnarray*}
   && {\cal{F}}=\{B_{r,x}:~~x\in A_1,~0<r<\delta,~B_{r,x}\subseteq\Omega,\\
   &&~~~~~~\mbox{and}~~\int_{B_{r,x}}|\nabla u|^2+e^u dx\ge \eta r\}
   \end{eqnarray*}
   By Vitali-Besicovitch Covering Theorem (see for instance \cite{EG}), we can find an at most contable family of points $(x_0^i)_{i\in I}$, $x_0^i\in A_1$ and
   $0<r_i<\delta$ such that $B_{r_i}(x_0^i)\in{\cal{F}} $ and $A_1\subseteq\cup_{i\in I}B_{r_i}(x_0^i)\,.$ Moreover every  $x\in A_1$ is contained in at most $N$
   balls, $N$ being a number depending only on the dimension of the space\,. \par
   The following estimates holds
   \begin{equation}\label{haus1}
  \eta \sum_{i\in I} r_i \le \sum_{i\in I} \int_{B_{r_i}(x_0^i)}|\nabla u|^2+e^u dx\,.
   \end{equation}
   \begin{eqnarray}\label{haus2}
 \sum_{i\in I}   \int_{B_{r_i}(x_0^i)}|\nabla u|^2+e^u dx&\le&\int_\Omega\sum_{i\in I}\11_{B_{r_i}(x_0^i)}(x)(|\nabla u|^2+e^u)dx\\
   &\le& N\int_{\{x:~d(x,A_1)\le \delta\}} (|\nabla u|^2+e^u) dx \le C\,,\nonumber
   \end{eqnarray}
   where $C>0$ is a constant independent on $\delta$.\par
   By combining (\ref{haus1}) and (\ref{haus2}) and letting $\delta\to 0$ we get that  ${\cal{H}}^1(A_1)=0$\,.
   \par
  {\bf  2.}  ${\cal{H}}^{1+\alpha}(A_2)=0\,.$\par
   Let $x\in A_2$. By definition 
   \begin{equation}\label{condbis}
   (u)_{r,x}{\cal{E}}_{r,x}(u)\ge\eta 
   \end{equation}
    for all $0<r<dist(x,\Omega)$\,.\par Jensen's Inequality implies that
    $$
    e^{(u)_{r,x}}\le \frac{1}{|\partial B_{r,x}|}\int_{\partial B_{r,x}} e^ u dx\,.
    $$
    Thus
    $$ (u)_{r,x}\le -C\log(r^2)\,.$$
    Thefore if $x\in A_2$ then for all $0<r<0$ we have $-C\log(r^2){\cal{E}}_{r,x}(u)\ge \eta\,.$\par
    Now fix $\delta \in (0,1)$. We have $A_2\subset\bigcup_{x\in A} B_{\delta ,x}\,.$ \par
     By Vitali's Covering Theorem
    there exists a countable number of disjoint balls $B_{5\delta,x_i}$, $x_i\in A_2,$ $i\in I$ such that
    $$A_2\subset\bigcup_{i\in I} B_{5\delta ,x_i}\,.$$
    We have
    \begin{eqnarray*}
    \int_\Omega(e^u+|\nabla u|^2) dx&\ge&\sum_{i\in I}\int_{B_{5\delta,x_i}}(e^u+|\nabla u|^2) dx\\
    &\ge& \eta C\sum_{i\in I}(5\delta)(\log{(5\delta)^{-2}})^{-1} \,.
    \end{eqnarray*}
    This implies that for all $0<\theta\le 1$ 
    $\sum_{i\in I}(5\delta)^\theta<+\infty$, hence by definition ${\cal{H}}_{dim}(A_2)\le 1\,.$
    \par
    It follows that ${\cal{H}}_{dim}(V)\le 1$ and we conclude\,.~~\hfill$\Box$ 
    \par
    \vskip1cm   
 {\bf Acknowledgements :} { The author is very grateful to Tristan Rivi\`ere for having drawn her attention to this problem and for having explained  
 the details of the technique he introduced together with Fang-Hua Lin in  \cite{LR1}. } 
 
 \par

\end{document}